\documentclass{article}
\usepackage{graphics}
\usepackage{latexsym}
\usepackage{amsfonts}
\usepackage{amssymb}
\usepackage{amsthm}
\usepackage{amsmath}
\usepackage{hyperref}
\author{Sheng Zhang}
\date{}
\title{Five-Collinear Sierpinski Gasket is\\Not Self-Similar}

\begin{document}
\newtheorem{theorem}{Theorem}[section]
\newtheorem{proposition}[theorem]{Proposition}
\newtheorem{corollary}[theorem]{Corollary}
\newtheorem{lemma}[theorem]{Lemma}
\maketitle

\begin{abstract}
A geometrical conclusion: Sierpinski gasket, two Sierpinski gaskets in a line, three Sierpinski gaskets in a line, and four Sierpinski gaskets in a line are self-similar, but five Sierpinski gaskets in a line is not, which is proved in this paper.
\end{abstract}

\section{Introduction}
Sierpinski gasket, two Sierpinski gaskets in a line, three Sierpinski gaskets in a line, and four Sierpinski gaskets in a line are the attractor of some contractive IFS consisting of similitudes, but five Sierpinski gaskets in a line is not. The proof is based on induction. Five Sierpinski gaskets in a line is of fractal dimension, which makes the situation a little bit complicated. The idea of the proof is to use figures of similar shapes with five Sierpinski gaskets in a line, but of integral dimension, to analyze properties of five Sierpinski gaskets in a line.

\section{Notations And Definitions}
\subsection{Notations}
\textbf{$\mathbb{Z}$} = $\{\dots, -2, -1, 0, 1, 2, \dots \}$.\\
\textbf{$\mathbb{N}$} = $\{0, 1, 2, \dots \}$.\\
\textbf{$\mathbb{N}_+$} = $\{1, 2, \dots \}$.\\
\textbf{$\blacktriangle P_i P_j P_k$}: The solid triangle in $\mathbb{R}^2$ with vertices $P_i$, $P_j$, and $P_k$.\\
\\
Below, let $A$, $B$ be sets and $f$, $g$ be maps:\\
\textbf{$d(x,y)$}: The distance between point $x$ and point $y$.\\
\textbf{$A-B$} = $\{x\mid x\in A, x\not\in B\}$.\\
\textbf{$diam(A)$} = $\sup_{x\in A, y\in A}\{d(x,y)\}$, where $A$ is a nonempty subset in $\mathbb{R}^2$.\\
\textbf{$d(A,B)$} = $\inf_{x\in A, y\in B}\{d(x,y)\}$, where $A$, $B$ are nonempty subsets in $\mathbb{R}^2$.\\
\textbf{$f\circ g$}: The composition of $f$ and $g$, which maps $x$ to $f(g(x))$.\\
\textbf{$f^{\circ k}$}: The composition of $k$ $f$'s ($k\in \mathbb{Z}$).

\subsection{Definitions}
\label{def}
\textbf{Similitude}: a map $\mathbb{R}^2\to \mathbb{R}^2$ which is the composition of scaling, rotation, translation, and maybe reflection. That is, $f$ is a similitude if and only if there exist $\theta \in [0,2\pi), k \in (0, +\infty)$ and $x_0, y_0\in \mathbb{R}$, such that
\begin{displaymath}
f{x \choose y}=k
\left( \begin{array}{cc}
\cos\theta & -\sin\theta \\
\sin\theta & \cos\theta
\end{array} \right)
{x \choose y} + {x_0 \choose y_0}
\end{displaymath}
or
\begin{displaymath}
f{x \choose y}=k
\left( \begin{array}{cc}
\cos\theta & \sin\theta \\
\sin\theta & -\cos\theta
\end{array} \right)
{x \choose y} + {x_0 \choose y_0},
\end{displaymath}
where $k$ is called the \textbf{scaling factor} of the similitude. If $k$ is strictly less than $1$, then the similitude is called $\textbf{contractive}$.\\
\textbf{IFS}(Iterated Function System)\cite{Hutchinson}: $F=\{ \mathbb{R}^2;f_1,f_2,\cdots,f_n\}$, where $f_1,f_2,\cdots,f_n:\mathbb{R}^2\to \mathbb{R}^2$ are continuous maps.\\
\textbf{Contractive IFS consisting of similitudes}: $F=\{ \mathbb{R}^2;f_1,f_2,\cdots,f_n\}$, where $f_1,f_2,\cdots,f_n:\mathbb{R}^2\to \mathbb{R}^2$ are contractive similitudes.\\
\textbf{Sierpinski gasket}: The attracor of the IFS $F=\{\mathbb{R}^2;
f_1{x \choose y}=
\left( \begin{array}{cc}
\frac{1}{2} & 0 \\
0 & \frac{1}{2}
\end{array} \right)
{x \choose y}, \\
f_2{x \choose y}=
\left( \begin{array}{cc}
\frac{1}{2} & 0 \\
0 & \frac{1}{2}
\end{array} \right)
{x \choose y}+{\frac{1}{4} \choose \frac{\sqrt{3}}{4}},
f_3{x \choose y}=
\left( \begin{array}{cc}
\frac{1}{2} & 0 \\
0 & \frac{1}{2}
\end{array} \right)
{x \choose y} + {\frac{1}{2} \choose 0}\}$ (see Figure \ref{figure1}).
\begin{center}
\makeatletter\def\@captype{figure}\makeatother
\scalebox{0.15}{\includegraphics{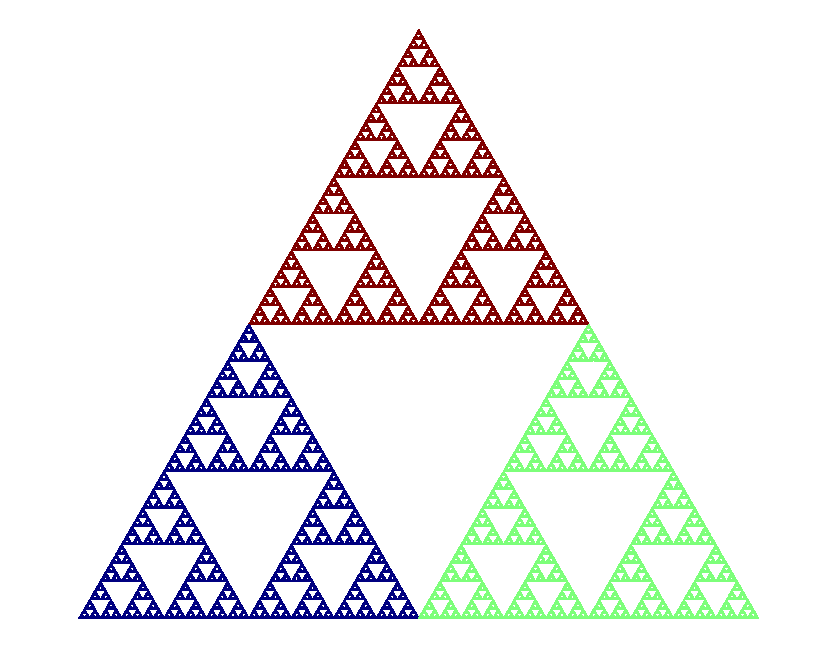}}
\caption{Sierpinski gasket}\label{figure1}
\end{center}
\textbf{Two-Sierpinski}: The union of Sierpinski gasket and Sierpinski gasket translated along positive $x$ axis by $1$ unit (see Figure \ref{figure2}).
\begin{center}
\makeatletter\def\@captype{figure}\makeatother
\scalebox{0.17}{\includegraphics{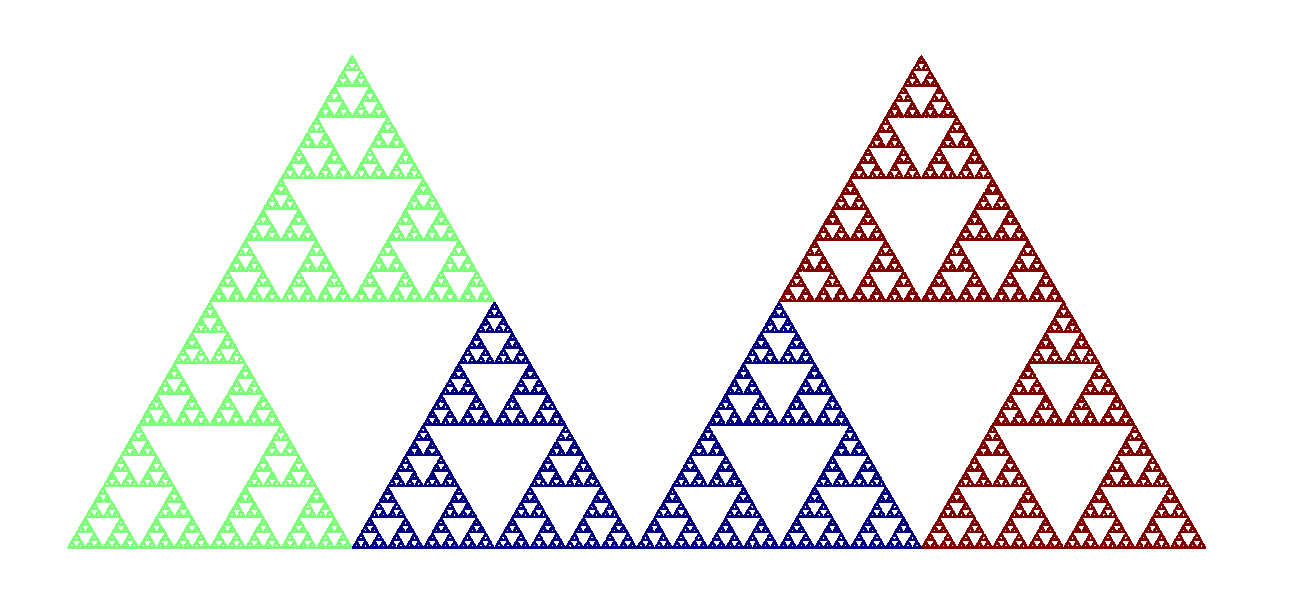}}
\caption{Two-Sierpinski}\label{figure2}
\end{center}
\textbf{$N$-Sierpinski}: The union of $(n-1)$-Sierpinski and Sierpinski gasket translated along positive $x$ axis by $n-1$ unit (see Figure \ref{figure3} when $n=3$). Five-Sierpinski will be denoted by \textbf{$E$} in this paper (see Figure \ref{figure4}). Let \textbf{$C$} be the figure obtained from Sierpinski gasket dilated by factor $8$. Then $E \subset C$ (see Figure \ref{figure5}, the whole figure is $C$, and the brown part is $E$).
\begin{center}
\makeatletter\def\@captype{figure}\makeatother
\scalebox{0.25}{\includegraphics{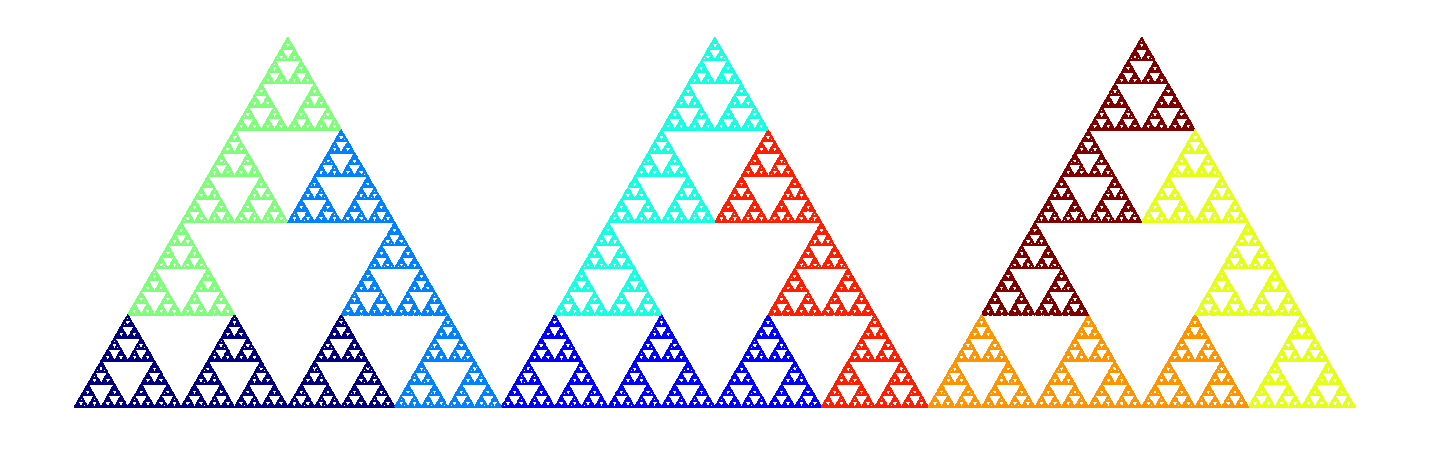}}
\caption{Three-Sierpinski}\label{figure3}
\end{center}
\begin{center}
\makeatletter\def\@captype{figure}\makeatother
\scalebox{0.27}{\includegraphics{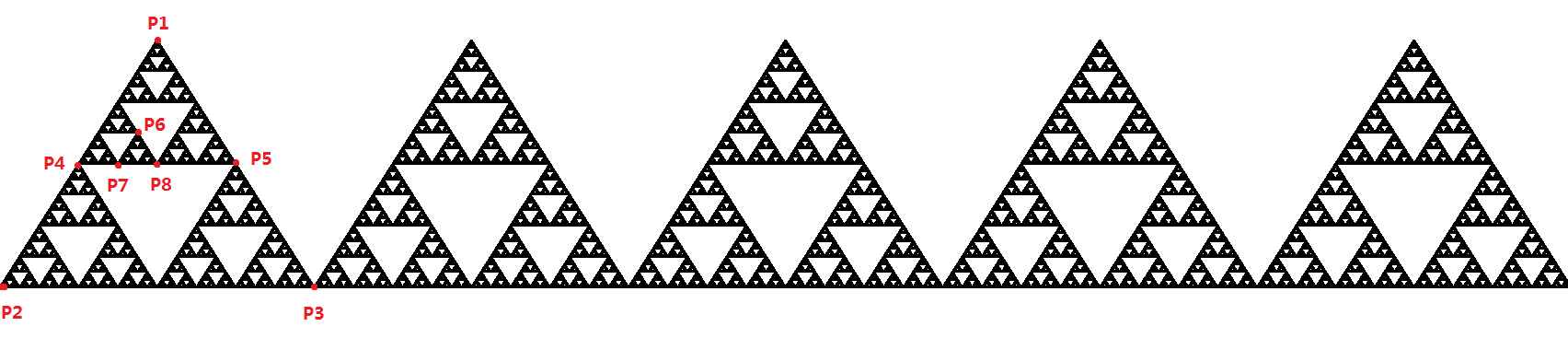}}
\caption{Five-Sierpinski (denoted by $E$)}\label{figure4}
\end{center}
\begin{center}
\makeatletter\def\@captype{figure}\makeatother
\scalebox{0.2}{\includegraphics{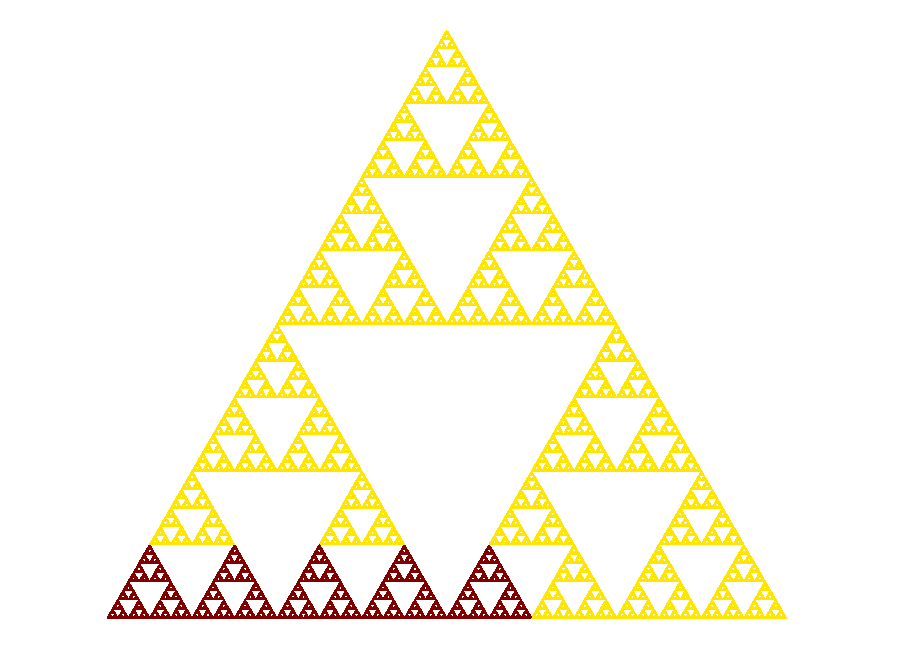}}
\caption{Five-Sierpinski (denoted by $E$) inside Sierpinski gasket dilated by factor $8$ (denoted by $C$)}\label{figure5}
\end{center}
\textbf{Figures $A_n$ and $B_n$} ($n\in \mathbb{Z}$ and $n \ge -3$): Since $C$ is a ``bigger'' version of Sierpinski gasket, it can be constructed as the intersection of a sequence of sets, the first four of which are listed in Figure \ref{figure6}. Denote this sequence of sets by $A_n$ ($n\in \mathbb{Z}$ and $n \ge -3$). Notice each $A_n$ is the union of $3^{n+3}$ solid equilateral triangles of the same size, whose topological interiors are disjoint each other. Denote the union of all vertices of these solid equilateral triangles by $B_n$ ($n\in \mathbb{Z}$ and $n \ge -3$), the first four of which are listed in Figure \ref{figure7}.
\begin{center}
\makeatletter\def\@captype{figure}\makeatother
\scalebox{0.7}{\includegraphics{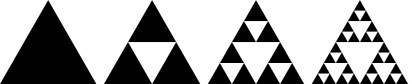}}
\caption{$A_{-3}$, $A_{-2}$, $A_{-1}$, and $A_{0}$}\label{figure6}
\end{center}
\begin{center}
\makeatletter\def\@captype{figure}\makeatother
\scalebox{0.2}{\includegraphics{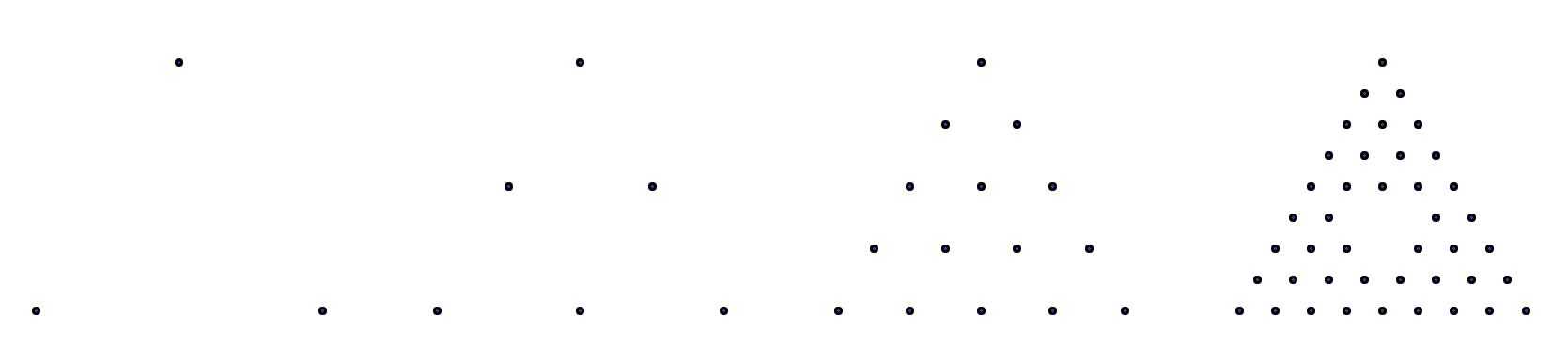}}
\caption{$B_{-3}$, $B_{-2}$, $B_{-1}$, and $B_{0}$}\label{figure7}
\end{center}
\textbf{Map $T$}: $\mathbb{R}^2 \to \mathbb{R}^2$,
\begin{equation}\label{def t}
T{x \choose y}=\frac{1}{2}
\left( \begin{array}{cc}
\cos\frac{2\pi}{3} & -\sin\frac{2\pi}{3} \\
\sin\frac{2\pi}{3} & \cos\frac{2\pi}{3}
\end{array} \right)
{x \choose y} + {\frac{3}{4} \choose \frac{\sqrt{3}}{4}}.
\end{equation}
Then $T$ is a contractive similitude with scaling factor $\frac{1}{2}$ and $T$ is bijective.

\section{Key Theorem}
As indicated in Figure \ref{figure1}, Figure \ref{figure2}, and Figure \ref{figure3}, Sierpinski gasket, two-Sierpinski, and three-Sierpinski are the attractor of some contractive IFS consisting of similitudes. These IFSs can be constructed as follows. In each figure, find similitudes that can map the whole figure to a region of the same colour. The number of similitudes is the number of different colours. Then the desired IFS is constructed by these similitudes. We can also construct a contractive IFS consisting of similitudes having four-Sierpinski as its attractor. This case is similar to three-Sierpinski.

However, as indicated in the following theorem, five-Sierpinski cannot be the attractor of any contractive IFS consisting of similitudes.
\begin{theorem}\label{theorem1}
Five-Sierpinski is not the attractor of any contractive IFS consisting of similitudes.
\end{theorem}
The proof of this theorem needs several lemmas.

\begin{lemma}\label{lemma1}
Denote five-Sierpinski by $E$. If $f$ is a similitude with scaling factor $k\le \frac{1}{80}$ such that $f(E)\subset E$ and $f(E)\cap \blacktriangle P_6 P_7 P_8 \neq \emptyset$ (see Figure \ref{figure4}), then $f(E)\subset \blacktriangle P_1 P_4 P_5 \cap E$. Further more, $T^{-1}(f(E))\subset \blacktriangle P_1 P_2 P_3 \cap E \subset E$.
\end{lemma}
\begin{proof}
First observe that $diam(E)=5$. Since $f$ is a similitude, for any $x,y\in \mathbb{R}^2$, $d(f(x),f(y))=k d(x,y)$. So, $diam(f(E))=\sup_{x\in f(E),y\in f(E)}\{d(x,y)\}=\sup_{x\in E,y\in E}\{d(f(x),f(y))\}=\sup_{x\in E,y\in E}\{k d(x,y)\}=k \sup_{x\in E,y\in E}\{d(x,y)\}=k\cdot diam(E)\le \frac{5}{80}=\frac{1}{16}$.

As $f(E)\cap \blacktriangle P_6 P_7 P_8 \neq \emptyset$, suppose $x_0 \in f(E)\cap \blacktriangle P_6 P_7 P_8$. Then, for all $x_1 \in f(E)$, $d(x_1,x_0)\le diam(f(E))=\frac{1}{16}$. Since $d(\blacktriangle P_6 P_7 P_8 , E-\blacktriangle P_1 P_4 P_5)=\frac{\sqrt{3}}{16}$, for all $x_2 \in (E-\blacktriangle P_1 P_4 P_5 )$, we have $d(x_2,x_0)\ge \frac{\sqrt{3}}{16}$. So, $f(E) \cap (E-\blacktriangle P_1 P_4 P_5 )= \emptyset$. Now, $f(E)\subset E$, so $f(E)\subset \blacktriangle P_1 P_4 P_5 \cap E$.

Since $\blacktriangle P_1 P_4 P_5 \cap E\subset T(E)$ ($T$ is defined in \ref{def t}), $f(E)\subset \blacktriangle P_1 P_4 P_5 \cap T(E)$. So, $T^{-1}(f(E))\subset T^{-1}(\blacktriangle P_1 P_4 P_5 \cap T(E))=T^{-1}(\blacktriangle P_1 P_4 P_5 )\cap E=\blacktriangle P_1 P_2 P_3 \cap E \subset E$.
\end{proof}
\begin{lemma}\label{lemma2}
Denote five-Sierpinski by $E$. If $f$ is a contractive similitude such that $f(E)\subset E$, then there exists $m\in \mathbb{N}_+$ such that the scaling factor of $f$ is $\frac{1}{2^m}$ and $f(P_1),f(P_2),f(P_3)\in B_m$ (See Figure \ref{figure4} and Figure \ref{figure7}).
\end{lemma}
\begin{proof}
As you can see in Figure \ref{figure5}, $E\subset C$, so $f(E)\subset E\subset C$. Because $C=\cap^{\infty}_{n=-3}A_n$, $f(E)\subset A_n$ for all $n\in \{ -3,-2,-1,0,\cdots \}$.

Each $A_n$ is the union of $3^{n+3}$ solid equilateral triangles of the same size, whose topological interiors are disjoint each other (see Figure \ref{figure6}). Consider three points $f(P_1),f(P_2),f(P_3)$. Now $f(P_1),f(P_2),f(P_3)$ are different points since $f$ is bijective. They are either in the same solid triangle or not. Suppose $N$ is the greatest integer in $\{ -3,-2,-1,0,\cdots \}$ such that $f(P_1),f(P_2),f(P_3)$ are in the same solid triangle of $A_N$ (note that such $N$ exists because $f(P_1),f(P_2),f(P_3)$ have positive distance each other and the diameter of solid triangles of $A_n$ tends to $0$ as $n$ tends to $\infty$). Then $f(P_1),f(P_2),f(P_3)$ are not in the same solid triangle of $A_{N+1}$.

Suppose $\blacktriangle Q_1 Q_4 Q_6$ is the solid triangle of $A_N$ where $f(P_1),f(P_2),f(P_3)$ are in (see Figure \ref{figure8}). Then $\blacktriangle Q_1 Q_2 Q_3$, $\blacktriangle Q_2 Q_4 Q_5$, $\blacktriangle Q_3 Q_5 Q_6$ are solid triangles of $A_{N+1}$. So $f(P_1),f(P_2),f(P_3)$ are not in the same one of $\blacktriangle Q_1 Q_2 Q_3$, $\blacktriangle Q_2 Q_4 Q_5$, $\blacktriangle Q_3 Q_5 Q_6$.
\begin{center}
\makeatletter\def\@captype{figure}\makeatother
\scalebox{1}{\includegraphics{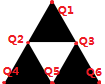}}
\caption{$\blacktriangle Q_1 Q_4 Q_6$}\label{figure8}
\end{center}

Suppose two of $f(P_1),f(P_2),f(P_3)$ are in the same one of $\blacktriangle Q_1 Q_2 Q_3$, $\blacktriangle Q_2 Q_4 Q_5$, $\blacktriangle Q_3 Q_5 Q_6$ and the other one of $f(P_1),f(P_2),f(P_3)$ is in another one of $\blacktriangle Q_1 Q_2 Q_3$, $\blacktriangle Q_2 Q_4 Q_5$, $\blacktriangle Q_3 Q_5 Q_6$. Without loss of generality, suppose $f(P_1),f(P_2)\in \blacktriangle Q_2 Q_4 Q_5$ and $f(P_3)\in \blacktriangle Q_3 Q_5 Q_6$. Since $f$ is a similitude, $f$ maps straight lines to straight lines. Because $f(E)\subset A_{N+1}$, $f$ maps segment $P_1 P_3$ to a segment in $A_{N+1}$ and maps segment $P_2 P_3$ to another segment in $A_{N+1}$. Now, segment $f(P_1)f(P_3)\subset A_{N+1}$ and segment $f(P_2)f(P_3)\subset A_{N+1}$.

If $f(P_3)=Q_3$, then $f(P_1)=Q_2$ and $f(P_2)=Q_5$, or $f(P_1)=Q_5$ and $f(P_2)=Q_2$. Both cases are impossible since $f(P_8)\in f(E)\subset A_{N+1}$ (see Figure \ref{figure4} for $P_8$). If $f(P_3)\neq Q_3$, then $f(P_1),f(P_2),f(P_3)$ are in segment $Q_4 Q_6$, which is also impossible.

Therefore, $f(P_1),f(P_2),f(P_3)$ are in different ones of $\blacktriangle Q_1 Q_2 Q_3$, $\blacktriangle Q_2 Q_4 Q_5$, $\blacktriangle Q_3 Q_5 Q_6$. Without loss of generality, suppose $f(P_1)\in \blacktriangle Q_1 Q_2 Q_3$, $f(P_2)\in \blacktriangle Q_2 Q_4 Q_5$, $f(P_3)\in \blacktriangle Q_3 Q_5 Q_6$. Then, using similar arguments as above, we can prove $f(P_1)=Q_1$, $f(P_2)=Q_4$, $f(P_3)=Q_6$.

Now, the scaling factor of $f$ is $\frac{1}{2^N}$ and $f(P_1),f(P_2),f(P_3)\in B_N$. As $f$ is contractive, $N\ge 1$. So, the lemma is proved.
\end{proof}
\begin{lemma}\label{lemma3}
For all $k\in \mathbb{N}_+$,
$$T^{\circ (k+1)}(\blacktriangle P_6 P_7 P_8)\subset T^{\circ k}(\blacktriangle P_6 P_7 P_8)\subset \cdots \subset \blacktriangle P_6 P_7 P_8$$ and
$$T^{\circ (k+1)}(\blacktriangle P_6 P_7 P_8 \cap E)\subset T^{\circ k}(\blacktriangle P_6 P_7 P_8 \cap E)\subset \cdots \subset \blacktriangle P_6 P_7 P_8 \cap E.$$
\end{lemma}
\begin{proof}
First observe that $T(\blacktriangle P_6 P_7 P_8)\subset \blacktriangle P_6 P_7 P_8$. For all $k\in \mathbb{N}_+$, by applying $T$, $T^{\circ 2}$, $\cdots$, $T^{\circ k}$ to both sides of this relation respectively, we have $T^{\circ 2}(\blacktriangle P_6 P_7 P_8)\subset T(\blacktriangle P_6 P_7 P_8)$, $T^{\circ 3}(\blacktriangle P_6 P_7 P_8)\subset T^{\circ 2}(\blacktriangle P_6 P_7 P_8)$, $\cdots$, $T^{\circ (k+1)}(\blacktriangle P_6 P_7 P_8)\subset T^{\circ k}(\blacktriangle P_6 P_7 P_8)$. Thus, $T^{\circ (k+1)}(\blacktriangle P_6 P_7 P_8)\subset T^{\circ k}(\blacktriangle P_6 P_7 P_8)\subset \cdots \subset \blacktriangle P_6 P_7 P_8$.

Similarly, since $T(\blacktriangle P_6 P_7 P_8 \cap E)\subset \blacktriangle P_6 P_7 P_8 \cap E$, we have $T^{\circ (k+1)}(\blacktriangle P_6 P_7 P_8 \cap E)\subset T^{\circ k}(\blacktriangle P_6 P_7 P_8 \cap E)\subset \cdots \subset \blacktriangle P_6 P_7 P_8 \cap E$.
\end{proof}
\begin{lemma}\label{lemma4}
If $f$ is a contractive similitude such that $f(E)\subset E$ and the scaling factor of $f$ is $\frac{1}{2^m}$ for some $m\in \mathbb{N}_+$, then $T^{\circ m}(\blacktriangle P_6 P_7 P_8)\cap f(E)=\emptyset$.
\end{lemma}
\begin{proof}
When $m\le 6$, according to Lemma \ref{lemma2}, $f(P_1),f(P_2),f(P_3)\in B_m$. Since $f$ is a similitude in $\mathbb{R}^2$, it is uniquely determined by the images of three noncollinear points $P_1,P_2,P_3$. As $f(P_1),f(P_2),f(P_3)\in B_m$ and $B_m$ consists of finitely many points, the possible choices of $f(P_1),f(P_2),f(P_3)$ are finite. So, the possible choices of $f$ are finite. Check directly among all possible choices of $f$ and find the statement $T^{\circ m}(\blacktriangle P_6 P_7 P_8)\cap f(E)=\emptyset$ always holds.

For any integer $M\ge 6$, suppose the lemma holds when $m=M$. When $m=M+1$, $f$ has scaling factor $\frac{1}{2^m}=\frac{1}{2^{M+1}}\le \frac{1}{2^7}\le \frac{1}{80}$.

If $f(E)\cap \blacktriangle P_6 P_7 P_8 =\emptyset$, then the lemma is proved since $T^{\circ m}(\blacktriangle P_6 P_7 P_8)\subset \blacktriangle P_6 P_7 P_8$.

If $f(E)\cap \blacktriangle P_6 P_7 P_8 \neq \emptyset$, then according to Lemma \ref{lemma1}, $T^{-1}(f(E))\subset E$. Let $g=T^{-1}\circ f$. Then $g$ is a contractive similitude such that $g(E)\subset E$ and the scaling factor of $g$ is $\frac{1}{2^{m-1}}=\frac{1}{2^M}$. By the hypothesis, $T^{\circ M}(\blacktriangle P_6 P_7 P_8)\cap g(E)=\emptyset$. Since $T$ is bijective, $T^{\circ m}(\blacktriangle P_6 P_7 P_8)\cap f(E)=T^{\circ (M+1)}(\blacktriangle P_6 P_7 P_8)\cap T(g(E))=T(T^{\circ M}(\blacktriangle P_6 P_7 P_8)\cap g(E))=T(\emptyset)=\emptyset$.

The lemma is proved by induction.
\end{proof}
\begin{lemma}\label{lemma4'}
If $f$ is a contractive similitude such that $f(E)\subset E$, then there exists $K\in \mathbb{N}_+$ such that $T^{\circ K}(\blacktriangle P_6 P_7 P_8)\cap f(E)=\emptyset$.
\end{lemma}
\begin{proof}
Suppose $f$ has scaling factor $k$. According to Lemma \ref{lemma2}, there exists $m\in \mathbb{N}_+$ such that $k=\frac{1}{2^m}$. According to Lemma \ref{lemma4}, $T^{\circ m}(\blacktriangle P_6 P_7 P_8)\cap f(E)=\emptyset$. Let $K=m$ and the lemma is proved.
\end{proof}
\begin{lemma}\label{lemma4''}
If $f$ is a contractive similitude such that $f(E)\subset E$, then there exists $K\in \mathbb{N}_+$ such that for all $k\ge K$, $T^{\circ k}(\blacktriangle P_6 P_7 P_8)\cap f(E)=\emptyset$.
\end{lemma}
\begin{proof}
This lemma is proved directly from Lemma \ref{lemma3} and Lemma \ref{lemma4'}.
\end{proof}
Now, let's prove the key theorem.
\begin{proof}[Proof of Theorem \ref{theorem1}]
Suppose there exists a contractive IFS $F=\{\mathbb{R}^2 ;f_1 ,f_2 ,\cdots ,f_n \}$ ($f_1 ,f_2 ,\cdots ,f_n$ are contractive similitudes) such that $E$ is the attractor of $F$. Then $f_1(E)\cup f_2(E)\cup \cdots \cup f_n(E)=E$. So, for all $m=1,2,\cdots ,n$, $f_m(E)\subset E$. According to Lemma \ref{lemma4''}, there exist $K_m \in \mathbb{N}_+$ such that for all $k\ge K_m$, $T^{\circ k}(\blacktriangle P_6 P_7 P_8)\cap f_m(E)=\emptyset$. Let $K=\max \{K_1,K_2,\cdots K_n\}$. Then for all $m=1,2,\cdots ,n$, $T^{\circ K}(\blacktriangle P_6 P_7 P_8)\cap f_m(E)=\emptyset$. So, $T^{\circ K}(\blacktriangle P_6 P_7 P_8)\cap (f_1(E)\cup f_2(E)\cup \cdots \cup f_n(E))=\emptyset$, or
\begin{equation}
T^{\circ K}(\blacktriangle P_6 P_7 P_8)\cap E=\emptyset.\label{theorem1 eq1}
\end{equation}

According to Lemma \ref{lemma3}, $T^{\circ K}(\blacktriangle P_6 P_7 P_8 \cap E)\subset \blacktriangle P_6 P_7 P_8 \cap E\subset E$. Since $\blacktriangle P_6 P_7 P_8 \cap E\neq \emptyset$, $T^{\circ K}(\blacktriangle P_6 P_7 P_8 \cap E)\neq \emptyset$. Thus, $T^{\circ K}(\blacktriangle P_6 P_7 P_8 \cap E) \cap E=T^{\circ K}(\blacktriangle P_6 P_7 P_8 \cap E)\neq \emptyset$. As $T^{\circ K}(\blacktriangle P_6 P_7 P_8 \cap E)\subset T^{\circ K}(\blacktriangle P_6 P_7 P_8)$, $T^{\circ K}(\blacktriangle P_6 P_7 P_8)\cap E\neq \emptyset$, which contradicts \ref{theorem1 eq1}.

The contradiction implies five-Sierpinski is not the attractor of any contractive IFS consisting of similitudes.
\end{proof}


\begin{thebibliography}{99}

\bibitem{Hutchinson}
Hutchinson, {\em Fractals and Self-Similarity}, Indiana University Journal of Mathematics 30: 713-747, 1981.

\end{thebibliography}
\end{document}